\documentclass{ifacconf}

\usepackage{natbib}            %
\usepackage{graphicx}          %
\RequirePackage[OT1]{fontenc}
\RequirePackage{amsmath,amssymb}
\numberwithin{equation}{section}
\newtheorem{lemma}{Lemma}[section]
\newtheorem{remark}[lemma]{Remark}
\newtheorem{theorem}[lemma]{Theorem}
\newtheorem{corollary}[lemma]{Corollary}

\newtheorem{proposition}[lemma]{Proposition}
\newtheorem{definition}[lemma]{Definition}
\newcommand{\B}[1]{{\bf #1}}

\newcommand{\R}[1]{{\rm #1}}
\newcommand{\mB}[1]{{\mathbb{#1}}}

\begin{document}

\begin{frontmatter}

\title{A statistical and computational theory for robust and sparse Kalman smoothing
} %

\author[First]{Aleksandr Aravkin}
\author[Second]{James V. Burke}
\author[Third]{Gianluigi Pillonetto}

\address[First]{Department of Earth and Ocean Sciences, University of British Columbia, Vancouver, Canada (e-mail: saravkin@eos.ubc.ca).}
\address[Second]{Department of Mathematics, University of Washington, Seattle, USA (e-mail: burke@math.washington.edu)}
\address[Third]{Department of Information Engineering, University of Padova, Padova, Italy (e-mail: giapi@dei.unipd.it)}


\begin{abstract}                          
Kalman smoothers reconstruct the state of a
dynamical system starting from noisy output samples. While the
classical estimator relies on quadratic penalization of process deviations 
and measurement errors, 
extensions that exploit Piecewise Linear Quadratic (PLQ)
penalties have been recently proposed in the literature. These new formulations 
include smoothers robust with respect to outliers in the data, and smoothers that keep better track of fast system
dynamics, e.g. jumps in the state values. In addition to $L_2$, well
known examples of PLQ penalties include the $L_1$, Huber and Vapnik losses. 
In this paper, we use a dual representation for PLQ penalties to build a statistical
modeling framework and a computational theory for Kalman smoothing.\\ 
We develop a statistical framework by
establishing conditions required to interpret PLQ penalties 
as negative logs of true probability densities. Then, we present a computational
framework, based on interior-point methods, that solves the
Kalman smoothing problem with PLQ penalties 
and maintains the linear complexity in the size of the time series,
just as in the $L_2$ case. %
The framework presented extends the computational efficiency of the Mayne-Fraser and 
Rauch-Tung-Striebel algorithms to a much broader non-smooth setting,
and includes many known robust and sparse smoothers as special cases.
\end{abstract}

\end{frontmatter}

\section{Introduction}

Consider the following discrete-time linear state-space model
\begin{equation}
\label{LinearGaussModel}
\begin{array}{rcl}
    x_1&=&x_0 + w_1
   \\
    x_k & = & G_k x_{k-1}  + w_k,   \qquad k=2,3,\ldots,N
    \\
    z_k & = & H_k x_k      + v_k,    \quad \qquad k=1,2,\ldots,N
\end{array}
\end{equation}
where $x_k \in \mB{R}^n$ is the state, $x_0$ is known, $z_k \in
\mB{R}^m$ contains noisy output samples, $G_k$ and $H_k$ are known
matrices. Further, $\{w_k\}$ and $\{v_k\}$ are mutually independent
zero-mean random variables with
covariances given by  $\{Q_k\}$ and $\{R_k\}$, respectively.\\
The classical fixed-interval Kalman smoothing problem is to obtain
 the (unconditional) minimum variance linear estimator of
the states $\{x_k\}_{k=1}^N$ as a function of $\{z_k\}_{k=1}^N$. It
is well known that the structure of this estimator
is related to the following optimization problem %
\begin{equation}
\label{KSNonlinObjective} \min_{\{x_k\}} %
\sum_{k=1}^N \|z_k - H_k x_k\|_{R_k^{-1}}^2 + \|x_k - G_k
x_{k-1}\|_{Q_k^{-1}}^2
\end{equation}
where $G_1$ denotes the identity matrix and $\|a\|^2_{\Sigma}:=a^\top
\Sigma a$  for every column vector $a$. When data become available,
the solution can be computed by the classical Kalman smoother with the
number of operations linear in $N$. This procedure also provides the
minimum variance estimate
of the states when all the system noises are assumed to be Gaussian.\\
In many circumstances, linear estimators  that rely on quadratic penalization of
model deviation, such as (\ref{KSNonlinObjective}), lead to unsatisfactory results. For
instance, they are not robust with respect to the presence of
outliers in the data \citep{Hub,Aravkin2011,Farahmand2011} and may
have difficulties in reconstructing fast system dynamics, e.g. jumps
in the state values \citep{Ohlsson2011}. In addition,
sparsity-promoting regularization is often used in order to extract
from a large measurement or parameter vector a small subset having
greatest impact on the predictive capability of the estimate for
future data. This sparsity principle permeates many well known
techniques in machine learning and signal processing, such as feature
selection, selective shrinkage, and compressed sensing
\citep{Hastie90,LARS2004,Donoho2006}. In many circumstances, when
smoothing is considered, it can be interpreted as a sparse non
Gaussian prior distribution on the noises entering the system. In these cases,
the estimator (\ref{KSNonlinObjective}) is often replaced by
\begin{equation} \label{prob2}
\sum_{k=1}^N V\left(z_k - H_k x_k;R_k\right)  +  J\left(x_k - G_k x_{k-1};Q_k \right)
\end{equation}
where, for example, $V$ can be the Huber or the Vapnik's
$\epsilon$-insensitive loss, used in support vector regression
\citep{Vapnik98,Evgeniou99}, while $J$ may be the  $\ell_1$-norm, as
in the LASSO procedure \citep{Lasso1996}.\\
The interpretation of problems such as (\ref{prob2}) in terms of Bayesian estimation
has been extensively studied in the statistical and machine learning literature
in recent years and
probabilistic approaches used in the analysis of estimation and learning algorithms
can be found e.g. in \citep{McKayARD,Tipping2001,Wipf_IEEE_TIT_2011}.
Non-Gaussian model errors and priors leading to a great %
variety of loss and penalty functions are also reviewed in
\citep{Wipf_ARD_NIPS_2006} using convex-type and integral-type
variational representations, with the latter being related to
Gaussian scale mixtures.
The fundamental novelty in this work is that, rather than taking this approach, we start
with a particular class of losses, called PLQ penalties, well known
from optimization literature \citep{RTRW}. We establish conditions
which allow these losses to be viewed as negative logs of true densities,
ensuring that $w_k$ and $v_k$ in (\ref{LinearGaussModel}) come from
true distributions. This in turn allows us to
interpret the solution to the problem (\ref{prob2}) as a MAP
estimator when the loss functions $V$ and $J$ come from this
subclass of PLQ penalties. We will show that this subclass includes
the four key examples, the $L_2$, $L_1$, Huber, and Vapnik
penalties.\\
The density characterization of PLQ penalties is achieved using a
dual representation, which also underlies the development of algorithms for fitting models of
the form (\ref{prob2}). In particular, in the second part of the
paper we derive the conditions, complimentary to those needed to set
up the statistical framework, that allow the development of new and
computationally efficient Kalman smoothers designed using 
non-smooth penalties on the process and measurement deviations. 
Amazingly, it turns out that the interior point method used in \citep{Aravkin2011}
generalizes perfectly to the entire class of PLQ densities under a
simple verifiable non-degeneracy condition. Hence, the solutions of
all the PLQ Kalman smoothers can be computed with a number of
operations that scales linearly in $N$, as in the quadratic case.
Such theoretical foundation generalizes the results recently
obtained in
\citep{Aravkin2011,AravkinIFAC,Farahmand2011,Ohlsson2011}, framing
them as particular cases of the framework presented here.\\
The paper is organized as follows. In Section \ref{PLQP} we
introduce the class of PLQ convex functions, and provide the
conditions under which they can be interpreted as negative logs 
 of corresponding densities. In Section \ref{InteriorPointKS} we
present a new PLQ Kalman smoother theorem
that generalizes the well known 
Mayne-Fraser two-filter and the Rauch-Tung-Striebel
algorithm \citep{Gelb} to nonsmooth formulations. This theorem is obtained 
by solving the Karush-Kuhn-Tucker (KKT) system for PLQ penalties 
using interior point methods, and exploiting the state space structure 
to obtain the solution in linear time. The necessary 
results and proofs supporting the main theorems appear in the Appendix.
We end the paper with a few concluding remarks.

\section{Piecewise Linear Quadratic Penalties and Densities}
\label{PLQP}

\subsection{Preliminaries}

We recall a few definitions from
convex analysis.

\begin{itemize}
\item \index{affine hull} (Affine hull) Define the affine hull of any set
$S$, denoted by $\R{aff}\; S$, as the smallest affine set that
contains $S$.
\item (Cone) For any set $S$, denote by $\R{cone}\; S$ the set $\{ts | s \in S, t
\in \mB{R}_+\}$.
\item \index{polar} (Polar Cone) For any cone $K \subset \mB{R}^m$,
the polar of $K$ is defined to be
\[
K^\circ := \{v  | \langle v, w \rangle \leq 0 \; \forall \; w \in
K\}.
\]
\item (Horizon cone). The (convex) Horizon cone $C^{\infty}$
is the set of `unbounded directions' for $C$, 
i.e. $d\in C^\infty$ if for any point $\bar w \in C$
we have $\{d | \bar w + \tau d \in \R{cl}\; C\;
\forall\; \tau \geq 0\}$. 

\end{itemize}

\subsection{PLQ densities}

We now introduce the PLQ penalties and densities that are the focus of this paper.
\begin{definition}
\index{penalties!piecewise linear-quadratic} (piecewise
linear-quadratic penalties) \citep{RTRW}. For a nonempty polyhedral
set $U \subset \mB{R}^m$ and a symmetric positive-semidefinite
matrix $M\in \mB{R}^{m\times m}$ (possibly $M =0$), the function
$\theta_{U, M}: \mB{R}^m \rightarrow \overline{\mB{R}}$ defined by

\begin{equation}
\label{PLQbasic}
\theta_{U, M}(w) := \sup_{u \in U}\left\{\langle u,w \rangle -
\frac{1}{2}\langle u, Mu \rangle\right\}
\end{equation}

\noindent is proper, convex, and piecewise linear-quadratic. When $M
=0$, it is piecewise linear; $\theta_{U, 0}=\sigma_U$, the support
function of $U$. The effective domain of $\theta_{U, M}$, denoted by
$\R{dom}(\theta_{U, M})$, is the set of $w \in \mB{R}^m$ where
$\theta_{U, M}(w) < \infty$, and is given by $(U^\infty \cap \R{Null}(M))^\circ$.
\vspace{-.55cm}
\begin{flushright}
$\blacksquare$
\end{flushright}\end{definition}

In order to capture the full class of penalties of interest, we
consider injective affine transformations into
$\mB{R}^m$ of the form $b + By$. The requirements on $B$ therefore
are $m \geq n$ and $\mathrm{Null}(B) = \{0\}$. The final technical
requirement we impose is that $b \in \R{dom}\; \theta_{U, M}$.
\begin{definition}(PLQ penalties with shifts and transforms)
\label{generalPLQ} Using (\ref{PLQbasic}), 
define $\rho: \mB{R}^n \rightarrow \mB{R}$
as $\theta_{U,M}(b + By)$:
\begin{equation}
\begin{array}{rcl}
\rho_{U, M, b, B}(y) &:=& 
\sup_{u \in U}
\left\{ \langle u,b + By \rangle - \frac{1}{2}\langle u, Mu
\rangle \right\} \;
\end{array}
\end{equation}
\vspace{-.55cm}
\begin{flushright}
$\blacksquare$
\end{flushright}\end{definition}

The following result characterizes the effective domain of $\rho$ (see Appendix for proof).

\begin{theorem}
\label{domainCharTheorem}  Let $\rho$ denote $\rho_{U, M, B, b}(y)$, and
$K$ denote $U^\infty \cap \mathrm{Null}(M)$. 
Suppose $U \subset \mB{R}^m$ is a
polyhedral set, $y \in \mB{R}^n$, $b \in K^\circ$, $M \in
\mB{R}^{m\times m}$ is positive semidefinite,
and $B\in \mB{R}^{m \times n}$ is injective.  
Then we have \((B^\R{T}K)^\circ \subset
\mathrm{dom}(\rho) \) and \((B^\R{T}(K\cap -K))^\perp =
\mathrm{aff}(\mathrm{dom}(\rho))\).
\vspace{-.55cm}
\begin{flushright}
$\blacksquare$
\end{flushright}\end{theorem}

Note that the functions $\rho$ are still piecewise linear-quadratic.
All of the important examples mentioned before can be
represented in this way, as shown below.

\begin{figure} \label{HuberVapnikFig}
  \begin{center} 
  \begin{tabular}{cc} 
  \includegraphics[scale=0.29]{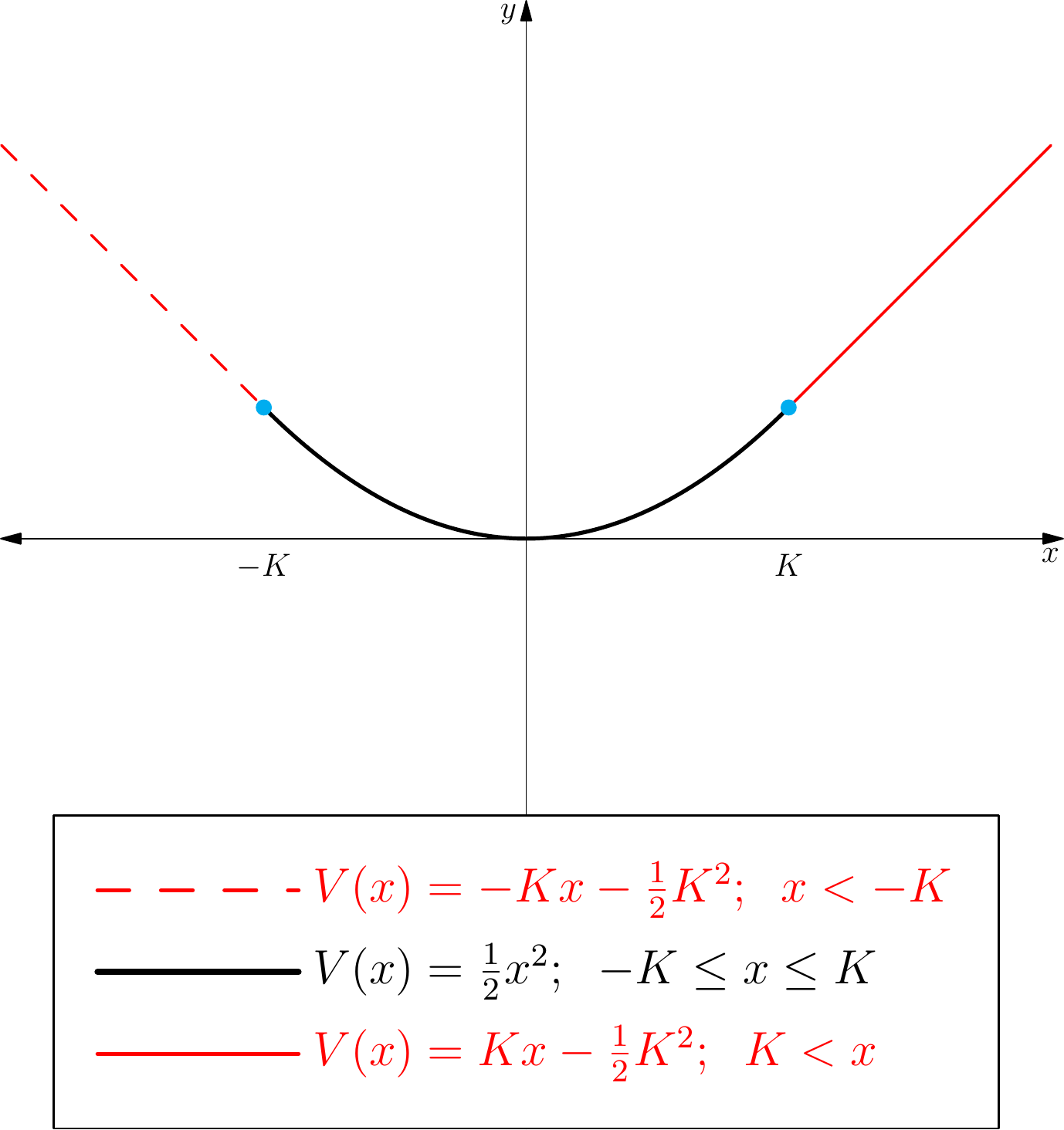}
\hspace{.1in}
{\includegraphics[scale=0.36]{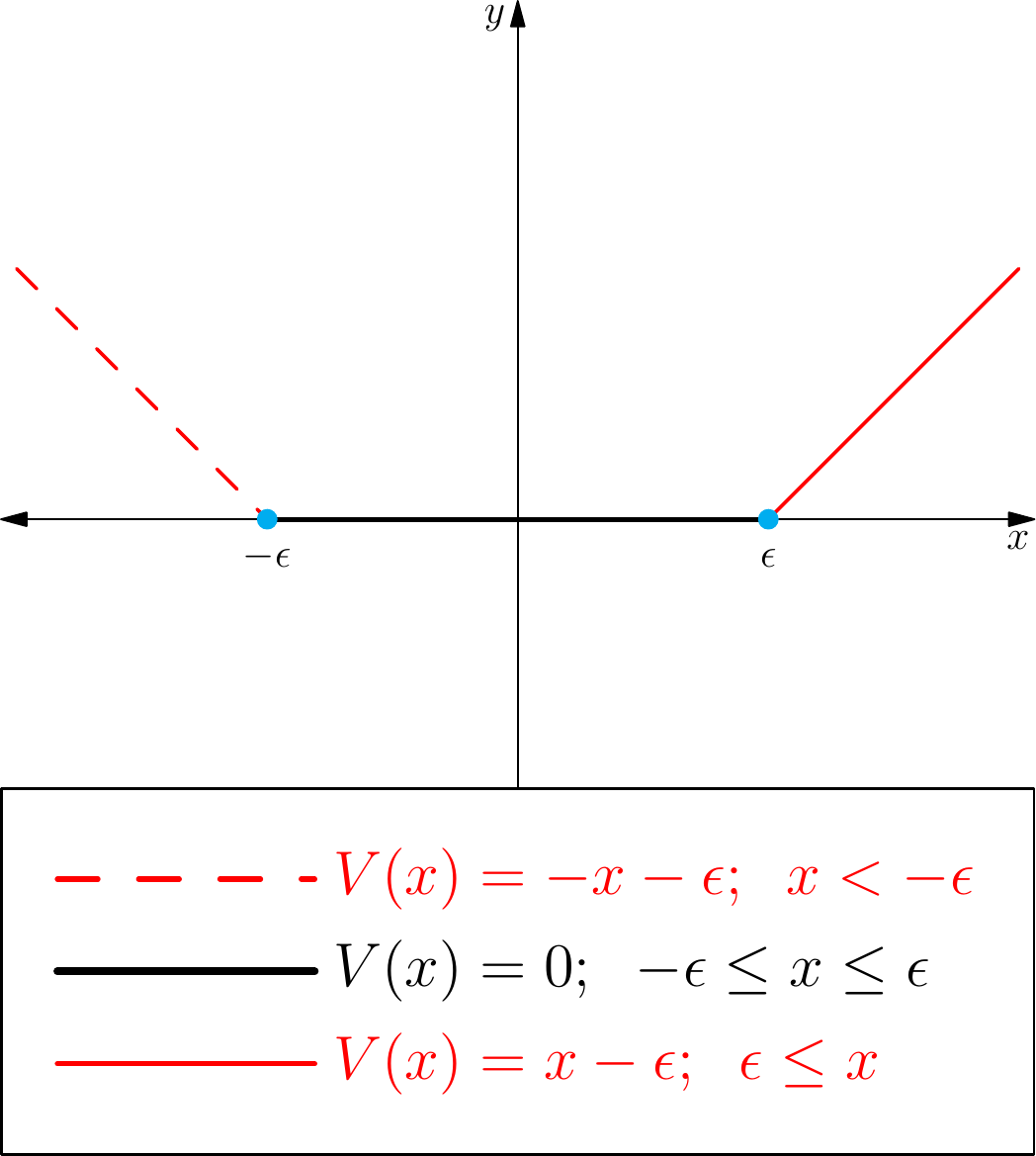}}
    \end{tabular}
    \caption{Huber (left) and Vapnik (right) Penalties}
     \end{center}
\end{figure}

\begin{remark}[scalar examples]
\label{scalarExamples} The $L_2$, $\ell_1$, Huber, and Vapnik
penalties are representable in the notation of Definition
\ref{generalPLQ}.
\begin{enumerate}
\item $L_2$: Take $U = \B{R}$, $M = 1$, $b = 0$, and $B = 1$. We obtain
\( \displaystyle \rho(y) = \sup_{u \in \B{R}}\left\langle uy -
\frac{1}{2}u^2 \right\rangle\;. \) The function inside the $\sup$ is
maximized at $u = y$, whence $\rho(y) = \frac{1}{2}y^2$.

\item $\ell_1$: Take $U = [-1, 1]$, $M = 0$, $b = 0$, and $B = 1$. We obtain
\( \displaystyle \rho(x) = \sup_{u \in [-1, 1]}\left\langle
uy\right\rangle\;. \) The function inside the $\sup$ is maximized by
taking $u = \R{sign}(y)$, whence $\rho(x) = |y|$.

\item Huber: Take $U = [-K, K]$, $M = 1$, $b = 0$, and $B = 1$.
We obtain \( \displaystyle \rho(y) = \sup_{u \in [-K,
K]}\left\langle uy - \frac{1}{2}u^2 \right\rangle\;. \) Take the
derivative with respect to $u$ and consider the following cases:
\begin{enumerate}
\item If $ y < -K $, take $u = -K$ to obtain
$-Ky  -\frac{1}{2}K^2$.
\item If $-K \leq y \leq K$, take $u = y$ to obtain
$\frac{1}{2}y^2$.
\item If $y > K $, take $u = K$ to obtain
a contribution of $Ky -\frac{1}{2}K^2$.
\end{enumerate}
This is the Huber penalty with parameter $K$, shown in the left
panel of Fig. 1.%

\item Vapnik: take $U = [0,1]\times[0,1]$,
$M = \left[\begin{smallmatrix}0 & 0\\0 & 0
\end{smallmatrix}\right]$, $B = \left[ \begin{smallmatrix} 1\\-1
\end{smallmatrix} \right]$, and $b = \left[ \begin{smallmatrix}
-\epsilon \\-\epsilon \end{smallmatrix} \right]$, for some $\epsilon
> 0$. We obtain \( \rho(y) = \sup_{u_1, u_2 \in [0,1]} \left\langle
\begin{bmatrix}
y-\epsilon\\
-y-\epsilon
\end{bmatrix},
\begin{bmatrix}
u_1\\
u_2
\end{bmatrix}
\right\rangle . \)
We can obtain an explicit representation by considering three cases:
\begin{enumerate}
\item If $|y| < \epsilon$, take $u_1 = u_2 = 0$. Then $\rho(y) = 0$.
\item If $y > \epsilon$, take $u_1 = 1$ and $u_2 = 0$. Then
$\rho(y) = y - \epsilon$.
\item If $y < -\epsilon$, take $u_1 = 0$ and $u_2 = 1$. Then
$\rho(y) = -y - \epsilon$.
\end{enumerate}
This is the Vapnik penalty with parameter $\epsilon$, shown in the
right panel of Fig. 1.%
\end{enumerate}
Note that the affine generalization (Definition \ref{generalPLQ}) is
already needed to express the Vapnik penalty.
\vspace{-.55cm}
\begin{flushright}
$\blacksquare$
\end{flushright}\end{remark}
In order to characterize PLQ penalties as negative logs 
of density functions, we need to ensure the integrability of said
density functions. A function $\rho(x)$ is called \emph{coercive} if
$\lim_{\|x\|\rightarrow \infty}\rho(x) = \infty$, and coercivity
turns out to be the key property to ensure integrability.
The proof of this fact, and the characterization of 
coercivity for PLQ penalties using the notation of 
Def. \ref{generalPLQ}, are the subject of the 
next two theorems (see Appendix for proofs).

\begin{theorem}
\label{PLQIntegrability} (PLQ Integrability). Suppose $\rho(y)$ is
coercive, and let $n_{\R{aff}}$ denote the dimension of
$\R{aff}(\R{dom}\; \rho)$. Then the function $f(y) = \exp(-\rho(y))$
is integrable on $\R{aff}(\R{dom}\; \rho)$ with the
$n_{\R{aff}}$-dimensional Lebesgue measure.
\vspace{-.55cm}
\begin{flushright}
$\blacksquare$
\end{flushright}\end{theorem}

\begin{theorem}
\label{coerciveRho} (Coercivity of $\rho$). $\rho$ is coercive if
and only if $[B^\R{T}\mathrm{cone}(U)]^\circ = \{0\}$.
\vspace{-.55cm}
\begin{flushright}
$\blacksquare$
\end{flushright}\end{theorem}

Theorem \ref{coerciveRho} can be used to show the coercivity
of familiar penalties. 
\begin{corollary} The penalties $L_2$, $L_1$, Vapnik, and
Huber are all coercive.
\end{corollary} %
{\bf{Proof:}} We show all of these penalties satisfy the hypothesis of Theorem
\ref{coerciveRho}.
\begin{enumerate}
\item $L_2$: $U = \B{R}$ and $B = 1$, so $[B^\R{T}\R{cone}(U)]^\circ =
\B{R}^\circ = \{0\}$.
\item $\ell_1$: $U = [-1, 1]$, so $\R{cone}(U) = \B{R}$, and $B = 1$,
so proof reduces to that case 1.
\item Huber: $U = [-K,K]$, so $\R{cone}(U) = \B{R}$, and $B = 1$, so proof reduces to that of case 1.
\item Vapnik: $U = [0,1] \times [0,1]$, so $\R{cone}(U) = \B{R}^2_+$.
$B = \left[ \begin{smallmatrix} 1\\-1 \end{smallmatrix} \right]$, so
$B^\R{T}\R{cone}(U) = \B{R}$, and again we reduce to case 1.
\end{enumerate}
\vspace{-.55cm}
\begin{flushright}
$\blacksquare$
\end{flushright}%
We now define a family of distributions on $\mB{R}^n$ by
interpreting piecewise linear quadratic functions $\rho$ as negative
logs of corresponding densities. 
Note that the support of the distributions is always contained 
in the affine set $\R{aff}(\R{dom}\; \rho)$, characterized 
in Th. \ref{domainCharTheorem}.

\begin{definition}
\label{PLQDensityDef} \index{density!piecewise linear-quadratic}
(Piecewise linear quadratic densities). Let $\rho$ be any coercive
piecewise linear quadratic function on $\mB{R}^n$ of the form
$\rho_{U, M, B, b;}(y) = \theta_{U, M}(b + By)$. Define $\B{p}(y)$
to be the following density on $\mB{R}^n$:
\begin{equation}
\label{PLQdensity} \B{p}(y) =
\begin{cases}
c_1^{-1}\exp(- \rho(y)) & y \in \R{dom}\; \rho\\
0 & \R{else},
\end{cases}
\end{equation}
where
\[
c_1 = \left(\int_{y \in \R{dom}\; \rho} \exp(-\rho(y))dy\right),
\]
and integral is with respect to the Lebesgue measure with
dimension $\R{dim}\Big(\R{aff}(\R{dom}\; \rho)\Big)$.
\vspace{-.55cm}
\begin{flushright}
$\blacksquare$
\end{flushright}\end{definition}

PLQ densities are true densities on the affine hull of the domain of $\rho$. %
The proof of Theorem \ref{PLQIntegrability} can be easily adapted to
show that they have moments of all orders. 
\section{Kalman Smoothing with PLQ penalties}
\label{InteriorPointKS}

In this section, we consider the model 
(\ref{LinearGaussModel}), but in the general case 
where errors $w_k$ and $v_k$ can come from any of the 
densities introduced in the previous section. To this 
end, we first formulate the KS problem over the entire 
sequence $\{x_k\}$.  

Given a sequence of column vectors $\{ u_k \}$ and matrices $ \{ T_k
\}$ we use the notation
\[
\R{vec} ( \{ u_k \} ) =
\begin{bmatrix}
u_1 \\ u_2  \\ \vdots \\ u_N
\end{bmatrix}
\; , \; \R{diag} ( \{ T_k \} ) =
\begin{bmatrix}
T_1    & 0      & \cdots & 0 \\
0      & T_2    & \ddots & \vdots \\
\vdots & \ddots & \ddots & 0 \\
0      & \cdots & 0      & T_N
\end{bmatrix} \; .
\]
We make the following definitions.
\[
\begin{array}{lll}
& x = \R{vec}\{x_1, \cdots, x_N\}\;,\qquad
&w = \R{vec}\{w_1, \cdots, w_K\}\; \\
&v = \R{vec}\{v_1, \cdots, v_k\}\;, \qquad
&Q = \R{diag}\{Q_1, \cdots, Q_N\}\;\\
&R = \R{diag}\{R_1, \cdots, R_N\}\;, \qquad
&H = \R{diag}\{H_1, \cdots, H_N\}.
\end{array}
\]

We also introduce the matrices $G$ and $H$:
\[
G =
\begin{bmatrix}
    \R{I}  & 0      &          &
    \\
    -G_2   & \R{I}  & \ddots   &
    \\
        & \ddots &  \ddots  & 0
    \\
        &        &   -G_N   & \R{I}
\end{bmatrix}\;,\; \quad
H =
\begin{bmatrix}
    H_1  & 0      &          &
    \\
    0   & H_2  & \ddots   &
    \\
        & \ddots &  \ddots  & 0
    \\
        &        &   0   & H_N
\end{bmatrix}\;.
\]

With this notation, model (\ref{LinearGaussModel}) can be written
\begin{equation}
\label{fullStat}
\begin{array}{lll}
\tilde x_0 
&=&
Gx + w\\
z 
&=&
Hx + v\;,
\end{array}
\end{equation}
where $x \in \mB{R}^{nN}$ is the entire state sequence of 
interest, $w$ is corresponding process noise, 
$z$ is the vector of all measurements, 
$v$ is the measurement noise, and
$\tilde x_0$ is a vector of size $nN$ with the first 
$n$-block equal to $x_0$, the initial state estimate, 
and the other blocks set to $0$. 

The general Kalman smoothing problem is 
described in the following proposition. 

\begin{proposition}
\label{prop:KS}
Suppose that the noises $w$ and $v$ in the model~\eqref{fullStat}
are PLQ densities with means $0$, variances  $Q$ and $R$
(see Def. \ref{PLQDensityDef}). Then, for suitable
$U^w, M^w,b^w,B^w$ and $U^v, M^v,b^v,B^v$ 
we have %
\begin{equation}
\label{kalmanDensities}
\begin{aligned}
\B{p}(w)&\propto \exp(-\theta_{U^w, M^w}(b^w + B^w Q^{-1/2}w)) \\
 \B{p}(v) &\propto \exp(-\theta_{U^v, M^v}(b^v + B^v R^{-1/2}v))\;
\end{aligned}
\end{equation}
while the MAP estimator
of $x$ in the model~\eqref{fullStat} is 
\begin{equation}
\arg \min_{x\in \mB{R}^{nN}} 
\left\{
\begin{aligned}
\label{PLQsubproblem} 
&\theta_{U^w,M^w}(b^w + B^w Q^{-1/2}(Gx - \tilde x_0)) \\
&+ \theta_{U^v, M^v}(b^v + B^vR^{-1/2}(Hx - z))
\end{aligned}
\right\}\; 
\end{equation} 
\vspace{-.55cm}
\begin{flushright}
$\blacksquare$
\end{flushright}\end{proposition}

Note that since $w_k$ and $v_k$ are independent,
problem~\eqref{PLQsubproblem} is decomposable into a sum of 
terms analogous to~\eqref{KSNonlinObjective}.
This special structure is manifest in the 
block diagonal structure of $H, Q, R, B^v, B^w$, 
the bidiagonal structure of $G$, and 
the structure of sets $U^w$ and $U^v$, 
and is key in proving the linear complexity
result that will be derived in the next part of this section.\\
For our purposes, it is now important to recall that, when the sets $U^w$ and $U^v$ are polyhedral, 
(\ref{PLQsubproblem}) is an Extended
Linear Quadratic program (ELQP), described in \cite[Example
11.43]{RTRW}. 
Rather than directly solving~\eqref{PLQsubproblem}, 
we work with the Karush-Kuhn-Tucker (KKT) system. 
We present the system in the following lemma, and derive
it in the Appendix. 
\begin{lemma}
\label{lem:KKT}
Suppose that the sets $U^w$ and $U^v$
are polyhedral, i.e. can be written 
\[
U^w = \{u|(A^w)^Tu \leq a^w \}, \quad U^v = \{u|(A^v)^Tu\leq a^v\}\;.
\]
Then the necessary first-order conditions 
for optimality of~\eqref{PLQsubproblem} are given by 
\begin{equation}
\label{PLQFinalConditions}
\begin{array}{lll}
&\begin{array}{llllll}
0 &=& (A^w)^\R{T}u^w + s^w - a^w\;;&&  0=  (A^v)^\R{T}u^v + s^v - a^v\\
0 &=& (s^w)^\R{T}q^w\;; && 0= (s^v)^\R{T}q^v
\end{array}\\\\
&\begin{array}{llllll}
0 &=& \tilde b^w + B^w Q^{-1/2}G\bar{d} -  M^w \bar{u}^w - A^wq^w\\
0 &=& \tilde b^v - B^v R^{-1/2}H\bar{d} - M^v \bar{u}^v -  A^v q^v\\
0 &=& G^\R{T}Q^{-\R{T}/2}(B^w)^\R{T}\bar u^w -
H^\R{T}R^{-\R{T}/2}(B^v)^\R{T}\bar u^v\\
0 &\leq& s^w, s^v, q^w, q^v. 
\end{array}
\end{array}
\end{equation}
\vspace{-.55cm}
\begin{flushright}
$\blacksquare$
\end{flushright}\end{lemma}

Our approach is to solve~\eqref{PLQFinalConditions}
directly  using Interior Point (IP) methods.  
IP methods work by applying a damped Newton iteration to a
relaxed version of (\ref{PLQFinalConditions}),  
specifically relaxing the `complementarity conditions':
\[
\begin{array}{lll}
(s^w)^\R{T}q^w = 0 & \rightarrow & Q^wS^w\B{1} - \mu\B{1} = 0 \\
(s^v)^\R{T}q^v = 0 & \rightarrow & Q^vS^v\B{1} - \mu\B{1} = 0\;,
\end{array}
\]
where $Q^w, S^w, Q^v, S^v$ are diagonal matrices
with diagonals $q^w, s^w, q^v, s^v$ respectively. 
The parameter $\mu$ is aggressively decreased to $0$ as the IP 
iterations proceed. Typically, no more than 10 or 20 iterations 
of the relaxed system are required to obtain a solution of~\eqref{PLQFinalConditions},
and hence an optimal solution to~\eqref{PLQsubproblem}.
The following theorem is key and represents the main result of this section. 
It shows that the computational effort required (per IP iteration)
is linear in the number of time steps whatever PLQ density 
enters the state space model. 

\begin{theorem}
\label{thm:PLQsmoother}
(PLQ Kalman Smoother Theorem) Suppose that all $w_k$ and $v_k$ in
the Kalman smoothing model (\ref{LinearGaussModel}) come from PLQ
densities that satisfy 
$\mathrm{Null}(M)\cap U^{\infty}  = \{0\}$, i.e. 
their corresponding penalties are finite-valued. 
Then~\eqref{PLQsubproblem} can be solved using an IP method, 
with computational complexity $O(Nn^3 + Nm)$ time.%
\vspace{-.55cm}
\begin{flushright}
$\blacksquare$
\end{flushright}\end{theorem}

The proof is presented in the Appendix 
and shows that IP methods for
solving~\eqref{PLQsubproblem} 
preserve the key block tridiagonal structure of the 
standard smoother. General smoothing
estimates can therefore be computed in $O(Nn^3)$ time, 
as long as the number of IP iterations is fixed 
(as it usually is in practice, to $10$ or $20$).\\
It is important to observe that the motivating examples 
(see Remark \ref{scalarExamples}) all satisfy the conditions 
of Theorem \ref{thm:PLQsmoother}.

\begin{corollary}
\label{cor:examples}
The densities corresponding to $L^1, L^2$, Huber, 
and Vapnik penalties 
all satisfy the hypotheses of Theorem \ref{thm:PLQsmoother}. 
\end{corollary}
{\bf Proof:}
We verify that $\mathrm{Null}(M) \cap \mathrm{Null}(A^\R{T}) = 0$
 for each of the four penalties. 
 In the $L^2$ case, $M$ has full rank.
 For the $L^1$, Huber, and Vapnik penalties, the respective
sets $U$ are bounded, so $U^{\infty}= \{0\}$.

\section{Conclusions}

We have presented a new theory for robust and sparse Kalman smoothing using
nonsmooth PLQ penalties applied to process and measurement deviations.
These smoothers can be designed within a statistical framework obtained 
by viewing PLQ penalties as  negative logs of true probability densities,
and we have presented necessary conditions that allow this interpretation. 
In this regard, the coercivity condition characterized in Th. \ref{coerciveRho}  
has been shown to play a central role. Notice that such a condition is also a nice example of how the
statistical framework established in the first part of this paper gives an alternative viewpoint for an idea
useful in machine learning. In fact, coercivity is also a
fundamental prerequisite in sparse and robust estimation as it precludes
directions for which the loss and the regularizer are insensitive to
large parameter/state changes. Thus, the condition for a (PLQ) penalty 
to be a negative log of a true density also
ensures that the problem is well posed and that the learning machine/smoother
can control model complexity.\\  
In the second part of the paper, we have 
shown that solutions to PLQ Kalman
smoothing formulations can be computed with a number of operations
that is linear in the length of the time series, as in the quadratic case.
A sufficient condition for the successful execution of IP iterations 
is that the PLQ penalties used should be finite valued, 
which implies non-degeneracy of the corresponding statistical
distribution (the support cannot be contained in a lower-dimensional
subspace). The statistical interpretation is thus strongly linked
to the computational procedure.\\ 
The computational framework presented allows 
a broad application of interior point methods to a wide class of 
smoothing problems of interest to
practitioners. The powerful algorithmic scheme designed here,
together with the breadth and significance of the new statistical framework presented,
underscores the practical utility and flexibility of this approach.
We believe that this perspective on model
development and Kalman smoothing will be useful in a number of applications
 in the years ahead.

\bibliographystyle{plain}
\bibliography{biblio}

\section*{Appendix}
\subsection*{Preliminaries}

\begin{definition}
\index{cone!horizon} (Horizon cone, specialized to the convex
setting by \cite[Theorem 3.6]{RTRW}). The Horizon cone $C^{\infty}$
for a convex set C is convex, and for any point $\bar w \in C$
consists of the vectors $\{d | \bar w + \tau d \in \R{cl}\; C\;
\forall\; \tau \geq 0\}$.
\end{definition}

\begin{definition}
\index{lineality} (Lineality). Define the lineality of convex cone
$K$, denoted $\R{lin}(K)$, to be $K \cap - K$. Since $K$ is a convex
cone, $\R{lin}(K)$ is the largest subspace contained in $K$.
\end{definition}

\begin{lemma}
\label{linealityLemma} (Characterization of lineality, \cite[Theorem
14.6]{RTR}). Let $K$ be any closed set containing the origin. Then
$\R{lin}(K) = (K^\circ)^{\perp}$.
\end{lemma}

\begin{definition}
\index{affine hull} (Affine hull). Define the affine hull of any set
$S$, denoted by $\R{aff}\; S$, as the smallest affine set that
contains $S$.
\end{definition}

\begin{corollary}
\label{affineHullCorollary} (Characterization of $\R{aff}\;
K^\circ$) Taking the perp of the characterization in Lemma
\ref{linealityLemma}, the affine hull of the polar of a closed
convex cone $K$ is given by $\R{aff} \; K^\circ = \R{lin}(K)^\perp$.
\end{corollary}

\subsection*{Proof of Theorem \ref{domainCharTheorem}}%

\begin{lemma}
\label{polarLemma} (Polars, linear transformations, and shifts) Let
$K\subset \mB{R}^n$ be a closed convex cone, $b \in \mB{R}^n$, and
$B \in \mB{R}^{n\times k}$. Then we have $(B^\R{T}K)^\circ \subset
B^{-1}(K^\circ - b)$ if $b \in K^\circ$.
\end{lemma}
{\bf{Proof:}} 
Recall that a convex cone is closed under addition. Then for any $b
\in K^\circ$, we have $K^\circ+b \subset K^\circ$, and hence
$K^\circ \subset K^\circ - b$. By \cite[Corollary 16.3.2]{RTR} we
get
\[
(B^\R{T}K)^\circ = B^{-1}K^\circ \subset B^{-1}(K^\circ - b)\;.
\]
\begin{flushright}
$\blacksquare$
\end{flushright}%

\begin{corollary}
\label{affineCorollary} Let $K$ be a closed convex cone, and $B\in
\mB{R}^{n\times k}$. If $b \in K^\circ$, then \(
\left(B^\R{T}(\mathrm{lin}(K))\right)^\perp \subset
\mathrm{aff}(B^{-1}(K^\circ - b)). \)
\end{corollary}
{\bf{Proof:}} %
By Lemma \ref{polarLemma}, \( \mathrm{aff}(B^{-1}(K^\circ - b))
\supset \mathrm{aff}(B^\R{T}K)^\circ =
\left(\mathrm{lin}(B^\R{T}K)\right)^\perp \) where the last equality
is by Corollary \ref{affineHullCorollary}. Since $B^\R{T}$ is a
linear transformation, we have $\mathrm{lin}(B^\R{T}K) =
B^\R{T}\mathrm{lin}(K)$.
\begin{flushright}
$\blacksquare$
\end{flushright}%

\begin{lemma}
\label{affineLemma} Let $K\subset \mB{R}^n$ be a closed convex cone,
$b \in \mathrm{aff}(K)^\circ$, and $B \in \mB{R}^{n\times k}$. Then
$\mathrm{aff}(B^{-1}(K^\circ - b)) \subset
B^{-1}(\mathrm{lin}(K))^\perp \subset (B^{-1}\mathrm{aff}(K^\circ -
b)) $.
\end{lemma}
{\bf{Proof:}} 
If $w \in \mathrm{aff}\left(B^{-1}(K^\circ - b)\right)$, for some
finite $N$ we can find sets $\{\lambda_i\}\subset \mB{R}$ and
$\{w_i\}\subset B^{-1}(K^\circ - b)$ such that $\sum_{i=1}^N
\lambda_i = 1$ and $\sum_{i=1}^N \lambda_iw_i = w$. For each $w_i$,
we have $Bw_i \in K^\circ - b$, so $b + Bw_i \in K^\circ$. Then
\[
b + Bw = \sum_{i = 1}^N \lambda_i (b + Bw_i) \in
\mathrm{aff}(K^\circ) = \mathrm{lin}(K)^\perp.
\]
Since $b\in \mathrm{lin}(K)^\perp$ by assumption, we have $Bw \in
\mathrm{lin}(K)^\perp$, and so $w \in
B^{-1}(\mathrm{lin}(K)^\perp)$.

Next, starting with $w \in B^{-1}(\mathrm{lin}(K)^\perp)$ we have
$Bw \in \mathrm{lin}(K)^\perp$ and so $b + Bw \in
\mathrm{lin}(K)^\perp$ since $\mathrm{lin}(K)^\perp$ is a subspace
and $b \in \mathrm{lin}(K)^\perp$. Then for some finite
 $\tilde N$ we can find
sets $\{\lambda_i\}\subset \mB{R}$ and $\{v_i\}\subset K^\circ$ such
that $\sum_{i=1}^{\tilde N} \lambda_i = 1$ and $\sum_{i=1}^{\tilde
N} \lambda_iv_i = b + Bw$. Subtracting $b$ from both sides, we have
$\sum_{i=1}^{\tilde N} \lambda_i(v_i - b) = Bw$, so in particular
$Bw \in \mathrm{aff}(K^\circ - b)$. Then $w \in
B^{-1}\mathrm{aff}(K^\circ - b)$.
\begin{flushright}
$\blacksquare$
\end{flushright}%

\begin{theorem}
\label{polarTheorem} Let $K\subset \mB{R}^n$ be a closed convex
cone, $b \in \mB{R}^n$, and $B \in \mB{R}^{n\times k}$. If $b \in
K^\circ$, then $(B^\R{T}\mathrm{lin}(K))^\perp =
\mathrm{aff}\left(B^{-1}(K^\circ - b)\right) =
B^{-1}(\mathrm{lin}(K)^\perp)$.
\end{theorem}
{\bf{Proof:}} 
From Corollary \ref{affineCorollary} and Lemma \ref{affineLemma}, we
immediately have
\[(B^\R{T}\mathrm{lin}(K))^\perp
\subset \mathrm{aff}\left(B^{-1}(K^\circ - b)\right) \subset
B^{-1}(\mathrm{lin}(K)^\perp).
\]
Note that for any subspace $S$, $S^\perp = S^\circ$. Then by
\cite[Corollary 16.3.2]{RTR}, $(B^\R{T}\mathrm{lin}(K))^\perp =
B^{-1}(\mathrm{lin}(K)^\perp)$.
\begin{flushright}
$\blacksquare$
\end{flushright}%

The proof of Theorem \ref{domainCharTheorem} now follows from Lemma \ref{polarLemma} and Theorem
\ref{polarTheorem}.
\subsection*{Proof of Theorem \ref{PLQIntegrability}}
Using the characterization of a piecewise quadratic function from
\cite[Definition 10.20]{RTRW}, the effective domain of
 $\rho(y)$ can be represented as the union of finitely many
polyhedral sets $U_i$, relative to each of which $\rho(y)$ is given
by an expression of the form $\frac{1}{2} \langle y, A_iy\rangle +
\langle a_i, y \rangle + \alpha_i$ for some scalar $\alpha_i \in
\mB{R}$, vector $a_i \in \mB{R}^n$ and symmetric positive
semidefinite matrix $A_i \in \mB{R}^{n\times n}$. Since $\rho(y)$ is
coercive, we claim that on each unbounded $U_i$ there must be some
constants $N_i$ and $\beta_i > 0$ so that for $\|y\| \geq N_i$ we
have $\rho(y) \geq \beta_i \|y\|$. Otherwise, we can find an index
set $J$ such that $\rho(y_j) \leq \beta_j \|y_j\|$, where $\beta_j
\downarrow 0$ and $\|y_j\|\uparrow \infty$. Without loss of
generality, suppose $\frac{y_j}{\|y_j\|}$converges to
 $\bar y \in U_i^\infty$, by \cite[Theorem 8.2]{RTR}.
By assumption, $\frac{\rho(y_j)}{\|y_j\|}\downarrow 0$, and we have
\[
\frac{\rho(y_j)}{\|y_j\|} = \|y_j\|\left\langle\frac{y_j}{\|y_j\|},
A_i\frac{y_j}{\|y_j\|} \right\rangle + \left\langle a_i,
\frac{y_j}{\|y_j\|}\right\rangle + \frac{\alpha_i}{\|y_j\|}.
\]
Taking the limit of both sides over $J$ we see that
$\|y_j\|\left\langle\frac{y_j}{\|y_j\|}, A_i\frac{y_j}{\|y_j\|}
\right\rangle$ must converge to a finite value. But this is only
possible if $\langle \bar y, A_i \bar y \rangle = 0$, so in
particular we must have $\bar y \in \mathrm{Null}(A_i)$. Note also
that $\langle a_i, \bar y \rangle \leq 0$, by taking the limit over
$J$ of
\[
\frac{\rho(y_j)}{\|y_j\|}\geq \left\langle a_i,
\frac{y_j}{\|y_j\|}\right\rangle + \frac{\alpha}{\|y_i\|},
\]
\vspace{.1cm}

\noindent so for any $x_0 \in U_i$ and $\lambda > 0$ we have $x_0 +
\lambda \bar y \in U_i$ since $\bar y \in U_i^\infty$ and
\[
\rho(x_0 + \lambda \bar y) \leq \rho(x_0) + \alpha_i,
\]
so in particular $\rho$ stays bounded as $\lambda \uparrow \infty$
and cannot be coercive.

The integrability of $f(y)$ is now clear. First note that $f(y)$ is
bounded below by $0$. Recall that the effective domain of $\rho$ can
be represented as the union of finitely many polyhedral sets $U_i$,
and for each unbounded such $U_i$ we have shown $f(y) \leq
\exp[-\beta_i\|y\|]$ off of some bounded subset of $U_i$. An
elementary application of the bounded convergence theorem shows that
$f$ must be integrable.
\subsection*{Proof of Theorem \ref{coerciveRho}}

First observe that $[B^{-1}(\mathrm{cone}(U)]^\circ =
[B^\R{T}\mathrm{cone}(U)]^\circ$ by \cite[Corollary 16.3.2]{RTR}.

Suppose that $\hat y \in B^{-1} ((\R{cone}\; U)^\circ)$, and $\hat y
\neq 0$. Then $B\hat y \in \mathrm{cone}(U)$, and $B\hat y \neq 0$
since $B$ is injective, and we have
\[
\begin{array}{lll}
\rho(\tau \hat y) &=& \sup_{u \in U} \langle b + \tau B \hat y,
u\rangle -
\frac{1}{2}u^\R{T}M u  \\
&=& \sup_{u \in U} \langle b , u\rangle - \frac{1}{2}u^\R{T}M u +
\tau  \langle B\hat y, u \rangle
\\
&\leq & \sup_{u \in U} \langle b , u\rangle -
\frac{1}{2}u^\R{T}M u \\
&\leq & \theta_{U, M}(b),
\end{array}
\]
so $\rho(\tau \hat y)$ stays bounded even as $\tau \rightarrow
\infty$, and so $\rho$ cannot be coercive.

Conversely, suppose that $\rho$ is not coercive. Then we can find a
sequence $\{y_k\}$ with $\|y_k\| > k$ and a constant $K$ so that
$\rho(y_k) \leq K$ for all $k > 0$. Without loss of generality, we
may assume that $\frac{y_k}{\|y_k\|}\rightarrow \bar y$.

Then by definition of $\rho$, we have for all $u \in U$
\[
\begin{array}{lll}
&\langle b + By_k, u \rangle - \frac{1}{2}u^\R{T}Mu \leq K\\
& \langle b + By_k, u \rangle \leq K + \frac{1}{2}u^\R{T} M u\\
& \langle \frac{b + By_k}{\|y_k\|}, u \rangle \leq \frac{K}{\|y_k\|}
+ \frac{1}{2\|y_k\|}u^\R{T} M u
\end{array}
\]
Note that $\bar y \neq 0$, so $B \bar y \neq 0$. When we take the
limit as $k \rightarrow \infty$, we get $\langle B\bar y, u \rangle
\leq 0$. From this inequality we see that $B \bar y \in (\R{cone}\;
U)^\circ$, and so $\bar y \in B^{-1}((\R{cone}\; U)^\circ)$.

\subsection*{Proof of Lemma \ref{lem:KKT}}
The Lagrangian for (\ref{PLQsubproblem})
for feasible $(x, u^w, u^v)$ is
\begin{equation}
\label{PLQLagrangian}
\begin{aligned}
\small L(x, u^w, u^v) &= \left\langle
\begin{bmatrix}
\tilde b^w
\\ \tilde b^v
\end{bmatrix},
\begin{bmatrix}
u^w
\\ u^v
\end{bmatrix}
\right\rangle - \frac{1}{2}
\begin{bmatrix}
u^w \\
u^v
\end{bmatrix}^\R{T}
\begin{bmatrix}
M^w & 0 \\
0 & M^v
\end{bmatrix}
\begin{bmatrix}
u^w \\
u^v
\end{bmatrix}\\
&- \left\langle
\begin{bmatrix}
u^w\\
u^v
\end{bmatrix}\;,
\begin{bmatrix}
- B^wQ^{-1/2}G \\
B^vR^{-1/2}H
\end{bmatrix}
x \right\rangle\;
\end{aligned}
\end{equation}
where $\tilde b^w = b^w - B^wQ^{-1/2}\tilde x_0$ and 
$\tilde b^v = b^v - B^vR^{-1/2}z$.
The associated optimality conditions for feasible $(x, u^w, u^v)$
are given by
\begin{equation}
\label{PLQOptimalityConditions}
\begin{array}{lll}
&G^\R{T}Q^{-\R{T}/2}(B^w)^\R{T}\bar u^w - H^\R{T}R^{-\R{T}/2}(B^v)^\R{T}\bar u^v  = 0\\
&\tilde b^w - M^w \bar{u}^w + B^w Q^{-1/2}G\bar{x} \in N_{U^w}(\bar{u}^w)\\
&\tilde b^v - M^v \bar{u}^v - B^v R^{-1/2}H\bar{x} \in
N_{U^v}(\bar{u}^v)\;,
\end{array}
\end{equation}
where $N_C(x)$ denotes the normal cone to the set $C$ at the point $x$ 
(see \cite{RTR} for details). 

Since $U^w$ and $U^v$ are polyhedral, we can derive
explicit representations of the normal cones $N_{U^w}(\bar u^w)$ and $N_{U^v}(\bar u^v)$.
For a polyhedral set $U \subset \mB{R}^m$
and any point $\bar{u} \in U$, the normal cone $N_U(\bar{u})$ is
polyhedral. Indeed, relative to any representation
\[
U = \{u|A^\R{T}u \leq a\}
\]
and the active index set \( I(\bar{u}) := \{i| \langle A_{i},
\bar{u}\rangle = a_i\} \), where $A_i$ denotes the $i$th column of
$A$, we have
\begin{equation}
\label{NormalRep}
N_U(\bar{u}) = \left\{ 
\begin{aligned}
q_1A_1 + \dots + q_mA_m\;  | \;q_i \geq 0\;
\R{for} \; & i\in I(\bar{u})\\
 q_i = 0 \; \R{for} & i \not\in I(\bar{u})
\end{aligned}
 \right\} .
\end{equation}
Using~\eqref{NormalRep}, 
 Then  we may rewrite the optimality conditions
(\ref{PLQOptimalityConditions}) more explicitly as 
\begin{equation}
\label{PLQExpandedConditions}
\begin{aligned}
&G^\R{T}Q^{-\R{T}/2}(B^w)^\R{T}\bar u^w -
H^\R{T}R^{-\R{T}/2}(B^v)^\R{T}\bar u^v  = 0\\
& \tilde b^w - M^w \bar{u}^w + B^w Q^{-1/2}G\bar{d} = A^wq^w\\
&\tilde b^v - M^v \bar{u}^v - B^v R^{-1/2}H\bar{d} = A^v q^v \\
&\{q^v \geq 0 | q_i^v = 0\; \R{for}\; i \not\in I(\bar u^v)\}\\
&\{q^w \geq 0 | q_i^w = 0\; \R{for}\; i \not\in I(\bar u^w)\}\;\\
\end{aligned}
\end{equation}

Define slack variables $s^w \geq 0$ and $s^v \geq 0$ as follows:
\begin{equation}
\label{slack}
\begin{array}{lll}
s^w &=&  a^w - (A^w)^\R{T}u^w\\
s^v &=&  a^v - (A^v)^\R{T}u^v .
\end{array}
\end{equation}
Note that we know the entries of $q_i^w$ and $q_i^v$ are zero if and
only if the corresponding slack variables $s_i^v$ and $s_i^w$ are
nonzero, respectively. Then we have $(q^w)^\R{T}s^w = (q^v)^\R{T}s^v
= 0$. These equations are known as the complementarity conditions.
Together, all of these equations give system~\eqref{PLQFinalConditions}.
\subsection{Proof of Theorem \ref{thm:PLQsmoother}}
IP methods apply a damped Newton iteration to 
find the solution of the  relaxed KKT system $F_{\mu} = 0$, 
where
\[
\small F_{\mu} \left(
\begin{matrix}
s^w\\
s^v\\
q^w\\
q^v\\
u^w\\
u^v\\
x
\end{matrix}
\right) =
\begin{bmatrix}
(A^w)^\R{T}u^w + s^w - a^w\\
(A^v)^\R{T}u^v + s^v - a^v\\
D(q^w)D(s^w)\B{1} - \mu\B{1}\\
D(q^v)D(s^v)\B{1} - \mu\B{1}\\
\tilde b^w + B^w Q^{-1/2}Gd -  M^w u^w - A^wq^w\\
\tilde b^v - B^v R^{-1/2}Hd - M^v u^v -  A^v q^v\\
G^\R{T}Q^{-\R{T}/2}(B^w)^\R{T} u^w -
H^\R{T}R^{-\R{T}/2}(B^v)^\R{T}\bar u^v
\end{bmatrix}.
\]
This entails solving the system
\begin{equation}
\label{NewtonSystem} \small F_{\mu}^{(1)} \left(
\begin{matrix}
s^w\\
s^v\\
q^w\\
q^v\\
u^w\\
u^v\\
d
\end{matrix}
\right)
\begin{bmatrix}
\Delta s^w\\
\Delta s^v\\
\Delta q^w\\
\Delta q^v\\
\Delta u^w\\
\Delta u^v\\
\Delta d
\end{bmatrix}
= -F_{\mu} \left(
\begin{matrix}
s^w\\
s^v\\
q^w\\
q^v\\
u^w\\
u^v\\
d
\end{matrix}
\right),
\end{equation}
where  the derivative matrix $F_{\mu}^{(1)}$ is given by
\begin{equation}
\label{Fprime} \tiny 
\begin{bmatrix}
I & 0 & 0 & 0 & (A^w)^\R{T} & 0 & 0\\
0 & I & 0 & 0 & 0 & (A^v)^\R{T} & 0\\
Q^w & 0 & S^w & 0 & 0 & 0 & 0\\
0 & Q^v & 0 & S^v & 0 & 0 & 0 \\
0 & 0 & - A^w & 0 & -M^w & 0 & B^wQ^{-1/2}G \\
0 & 0 & 0 & -A^v & 0 & -M^v & -B^vR^{-1/2}H \\
0 & 0 & 0 & 0 & G^\R{T}Q^{-\R{T}/2}(B^w)^\R{T} &
-H^\R{T}R^{-\R{T}/2}(B^v)^\R{T} & 0
\end{bmatrix}
\end{equation}
We now show the row operations necessary to reduce the matrix
$F_{\mu}^{(1)}$ in (\ref{Fprime}) 
to upper block triangular form. After each
operation, we show only the row that was modified.
\begin{equation*}
\begin{array}{lll}
&\R{row}_3 \leftarrow \R{row}_3 - D(q^w)\;\R{row}_1\\
&\begin{bmatrix}
0 & 0 & D(s^w) & 0 & -D(q^w)(A^w)^\R{T} & 0 & 0\\
\end{bmatrix}\\
&\R{row}_4 \leftarrow \R{row}_4 - D(q^v)\;\R{row}_2 \\
&\begin{bmatrix}
0 & 0 & 0 & D(s^v) & 0 & -D(q^v)(A^v)^\R{T} & 0 \\
\end{bmatrix}\\
&\R{row}_5 \leftarrow \R{row}_5 + A^wD(s^w)^{-1}\;\R{row}_3\\
&\begin{bmatrix}
0 & 0 & 0 & 0 & -T^w & 0 & B^wQ^{-1/2}G \\
\end{bmatrix}\\
&\R{row}_6 \leftarrow \R{row}_6 + A^vD(s^v)^{-1}\;\R{row}_4 \\
&\begin{bmatrix}
0 & 0 & 0 & 0 & 0 & -T^v  & -B^vR^{-1/2}H
\end{bmatrix}\;.
\end{array}
\end{equation*}
In the above expressions,
\begin{equation}
\label{KalmanT}
\begin{aligned}
T^w &:=  M^w + A^w(S^w)^{-1}Q^w(A^w)^\R{T}\\
T^v &:= M^v + A^v(S^v)^{-1}Q^v(A^v)^\R{T}\;,
\end{aligned}
\end{equation}
where $(S^w)^{-1}Q^w$ and $(S^v)^{-1}Q^v$ are always 
full-rank diagonal matrices, since the vectors $s^w, q^w, s^v, q^v$
are always strictly positive in IP iterations. 
The invertibility of $T^w$ and $T^v$ is charachterized 
in the following lemma. 
\begin{lemma}
\label{TInvertibility} (Invertibility of $T$) Let $\theta_{U,
M}(\cdot)$ be any PLQ penalty on $\mB{R}^k$, with $U = \{u \Big| A^\R{T}u \leq
a\}$. Let $T_\theta := M + ADA^\R{T}$, where $D$ is any diagonal
$k\times k$ matrix with positive entries on the diagonal.
Then $T_\theta$ is invertible if and only if $\mathrm{Null}(M)\cap
U^{\infty}  = \{0\}$, or $\R{dom}(\theta_{U, M})$ is $\mB{R}^k$. 
\end{lemma}
{\bf Proof:}
Note that
\[
\begin{array}{lll}
\mathrm{Null}(M + ADA^\R{T}) &=&
\{w \Big| w^\R{T}Mw + w^\R{T}ADA^\R{T}w = 0\}\\
&=&
\{w \Big| w \in \mathrm{Null}(M)\;,\; w\in \mathrm{Null}(A^\R{T})\}\\
&=& \mathrm{Null}(M) \cap \mathrm{Null}(A^\R{T}).
\end{array}
\]
The first claim now follows from the fact that $U^{\infty} =
\mathrm{Null}(A^\R{T})$. 
Recall that the effective domain of $\theta$ is given by 
$(\R{Null}(M)\cap U^\infty)^\circ$, and it is immediate 
from the definition of `polar' that $0^\circ = \mB{R}^k$.
\vspace{-.55cm}
\begin{flushright}
$\blacksquare$
\end{flushright}\begin{remark}
\label{TblockDiagonal} (Block diagonal structure of $T$ in i.d.
case) Suppose that $\B{y}$ is a random vector, $\B{y} =
\left(\begin{smallmatrix} \B{y_1} & \cdots &
\B{y_n}\end{smallmatrix}\right)$, where each $\B{y_i}$ is itself a
random vector in $\mB{R}^{m_i}$, from some PLQ density $\B{p}(y_i)
\propto \exp[-c_2\theta_{U_i, M_i}((\cdot))]$, and all $\B{y_i}$
are independent. Let $U_i = \{u: A_i^Tu \leq a_i\}$. Then the matrix
$T_\theta$ is given by 
\(
T_\theta = M + ADA^T
\)
where 
$M = \R{diag}[M_1, \cdots, M_N]$, 
$A = \R{diag}[A_1, \cdots,A_N]$,
$D = \R{diag}[D_1, \cdots,D_N]$, 
and $\{D_i\}$ are diagonal with positive entries.
Moreover, $T_\theta$ is block diagonal, with $i$th diagonal block given by
$M_i + A_iD_iA_i^\R{T}$.
\vspace{-.55cm}
\begin{flushright}
$\blacksquare$
\end{flushright}
\end{remark}

\begin{corollary}
($T$ matrices in the Kalman smoothing context) The matrices $T^w$
and $T^v$ in (\ref{KalmanT}) are block diagonal provided that
$\{w_k\}$ and $\{v_k\}$ are independent vectors from any PLQ
densities. Moreover, if these densities all satisfy the
characterization in Lemma \ref{TInvertibility}, these matrices are
also invertible.
\vspace{-.55cm}
\begin{flushright}
$\blacksquare$
\end{flushright} 
\end{corollary}

We now finish the reduction of $F_{\mu}^{(1)}$ to upper block
triangular form:
\begin{equation*}
\begin{aligned}
 \R{row}_7 &\leftarrow \R{row}_7 +
\Big(G^\R{T}Q^{-\R{T}/2}(B^w)^\R{T}(T^w)^{-1}\Big)\R{row}_5 -\\
&\Big(H^\R{T}R^{-\R{T}/2}(B^v)^\R{T}(T^v)^{-1}\Big)\R{row}_6 \\
&\small\begin{bmatrix}
I & 0 & 0 & 0 & (A^w)^\R{T} & 0 & 0\\
0 & I & 0 & 0 & 0 & (A^v)^\R{T} & 0\\
0 & 0 & S^w & 0 & -Q^w(A^w)^\R{T} & 0 & 0\\
0 & 0 & 0 & S^v & 0 & -Q^v(A^v)^\R{T} & 0 \\
0 & 0 & 0 & 0 & -T^w & 0 & B^wQ^{-1/2}G \\
0 & 0 & 0 & 0 & 0 & -T^v  & -B^vR^{-1/2}H \\
0 & 0 & 0 & 0 & 0 & 0 & \Phi
\end{bmatrix}
\end{aligned}
\end{equation*}
where
\begin{equation}
\label{Phi} 
\begin{aligned}
\Phi &= \Phi_G + \Phi_H = 
G^\R{T}Q^{-\R{T}/2}(B^w)^\R{T}(T^w)^{-1}B^wQ^{-1/2}G \\
&+
H^\R{T}R^{-\R{T}/2}(B^v)^\R{T}(T^v)^{-1}B^vR^{-1/2}H.
\end{aligned}
\end{equation}
Note that $\Phi$ is symmetric positive definite. Note also that
$\Phi$ is block tridiagonal, since
\begin{enumerate}
\item $\Phi_H$ is block diagonal.
\item $Q^{-\R{T}/2}(B^w)^\R{T}(T^w)^{-1}B^wQ^{-1/2}$ is block
diagonal, and $G$ is block bidiagonal, hence $\Phi_G$
is block tridiagonal.
\end{enumerate}
Solving  system (\ref{NewtonSystem}) requires inverting the block
diagonal matrices $T^v$ and $T^w$ at each iteration of the damped
Newton's method, as well as solving an equation of the form $\Phi
\Delta x = \varrho$. We have already seen that $\Phi$ is block
tridiagonal, symmetric, and positive definite, so $\Phi
\Delta x = \varrho$ can be solved in $O(Nn^3)$ time using the block
tridiagonal algorithm in \citep{Bell2000}. The remaining four back
solves required to solve (\ref{NewtonSystem}) can each be done in
$O(nN)$ time.

\end{document}